\documentclass[a4paper,12pt]{amsart}

\usepackage{amsmath,amsfonts,amssymb,mathrsfs}

\DeclareMathOperator{\esp}{{\mathbb E}} 
\DeclareMathOperator{\proba}{{\mathbb P}} 

\newcommand*{\e}{\mathrm{e}}
\newcommand*{\dif}{{\mathsf d}}

\begin{document}

\title[Explicit fixed points of the smoothing tranformation]{Explicit
  fixed points of the smoothing tranformation}
\author{Jacques Peyri\`ere}

\address{Institut de Mathématiques d’Orsay, CNRS, Université
  Paris-Saclay, 91405 Orsay, France and Renmin University of China,
  School of Mathematics, 100872 Beijing, P.R.\ China.}

\email{peyriere@phare.normalesup.org}

 \date{}

\maketitle
\begin{abstract} We deal with the equation
  $\displaystyle Y \stackrel{\dif}{=} \frac{1}{b} \sum_{1\le j\le N}
  W_jY_j$, where the unknown is the distribution of~$Y$, the variables
  in the right hand side are independent, the $Y_j$ are
  equidistributed with $Y$, $N$ is an integer valued random variable,
  and the $W_j$ are equidistributed, nonnegative and of
  expectation~1. Usually a solution is obtained as the limit of a
  martingale. In some cases we give an explicit formula for the law
  of~$Y$.
\end{abstract}
  
Let $N$ be a nonnegative integer valued random variable with finite
second moment. Its expectation~$b$ is assumed to be larger than~1.
Set $\varphi(x) = \esp x^N = \sum_{n\ge 0} \proba(N=n)x^n$. Then
$b = \varphi'(1)$. Let~$\mu$ be a probability measure on $[0,+\infty)$ such that
\begin{equation}
  \int_0^\infty x\,\mu(\dif x) = 1 \text{\quad and\quad} \int_0^\infty x^2\mu(\dif x) < b.
\end{equation}

We consider the equation
\begin{equation}\label{equation}
  Y \stackrel{\dif}{=} \frac{1}{b} \sum_{1\le j\le N}
  W_jY_j
\end{equation}
whose unknown is the probability law~$\nu$ on $[0,+\infty)$ common
to~$Y$ and all the $Y_j$, the variables on the right hand side are
independent and all the variables~$W_j$ are distributed
according to~$\mu$.

 \subsection*{The Mandelbrot martingales}

 To get a solution to Equation~\eqref{equation} one way is to use the
 Mandelbrot construction~\cite{M1}. We set
\begin{equation*}
  Y_n = b^{-n}\sum_{1\le j_1\le N} W_{j_1}
\sum_{\substack{1\le j_i\le N_{j_1,j_2\dots,j_{i-1}}\\ \text{for }1\le
    i \le n}} W_{j_1,j_2}W_{j_1,j_2,j_3}\dots W_{j_1,j_2,\dots,j_n},
\end{equation*}
where all the variables in the right hand side are independent, the
$W_{j_1,j_2\cdots j_n}$ are distributed according to~$\mu$, and~$N$
and the variables $N_{j_1,j_2\cdots j_n}$ are equidistributed.

One has
\begin{multline*}
  Y_{n+1} = b^{-n}\sum_{1\le j_1\le N} W_{j_1} \times\\
  \sum_{\substack{1\le j_i\le N_{j_1,j_2\dots,j_{i-1}}\\ \text{for
      }1\le i \le n}} W_{j_1,j_2}W_{j_1,j_2,j_3}\dots
  W_{j_1,j_2,\dots,j_n}\ b^{-1}\hspace{-1.5em}\sum_{1\le j_{n+1}\le
    N_{n+1}} W_{j_1,j_2,\dots,j_{n+1}}.
\end{multline*}
Therefore $(Y_n)_{n\ge 1}$ is a martingale.

We also have
\begin{equation}
  Y_{n+1} =b^{-1} \sum_{1\le j_1\le N_1} W_{j_1}Y_{n}(j_1) \label{recur}
\end{equation}
where all the variables are independent and the $Y_n(j)$ are
equidistributed with~$Y_n$.

\noindent So,
$\esp Y_{n+1}^2 = b^{-2}\esp W^2 \esp Y_n^2 \esp N + b^{-2} \esp
N(N-1)$, i.e.,\\[3pt]
$ \esp Y_{n+1}^2 = b^{-1}\esp W^2 \esp Y_n^2 + b^{-2} \esp N(N-1)$.

We therefore see that if $\esp W^2< b$ the martingale $(Y_n)$ is
bounded in~$L^2$. Then it has a limit $Y$, almost surely and in $L^2$,
$\esp Y =1$, and
$$\esp Y^2 = \displaystyle \frac{\esp N(N-1)}{b(b-\esp W^2)}.$$

Due to~\eqref{recur} we see that~$Y$ is a solution to
Equation~\eqref{equation}.\medskip

B.~Mandelbrot~\cite{M1,M2} introduced this construction (with
constant~$N$) to give a simplified statistical description of the
dissipation of energy in a turbulent flow. Since then this model and
Equation~\eqref{equation} have been extensively studied and
generalized (for instance~\cite{KP,DL,G,L}). See~\cite{BFP,BP} for
a survey.\medskip

As a matter of fact, the necessary and sufficient condition for the
uniform integrability of the martingale~$Y_n$ is
$\esp W\log W< \log b$ (see~\cite{KP} when~$N$ is constant and
\cite{L} when~$N$ is not constant).\medskip

Let $\alpha = \proba(W\ne 0)$ and $\beta = \proba(Y=0)$. It results
from Equation~\eqref{equation} that
$\beta = \varphi(\alpha\beta + 1-\alpha)$.  Due to convexity, the
function $t\mapsto \varphi(\alpha t+1-\alpha)$ has a fixed point less
than~1 if and only if its derivative at~1 is larger than~1, which
means $\alpha b>1$. Due to hypotheses $\esp W = 1$ and $\esp W^2< b$
this condition is fulfilled.

\subsection*{Laplace transform}

For $t\ge 0$ set $g(t) = \esp \e^{-tY}$. Then
\begin{equation}\label{init}
g(0)=-g'(0)=1.
\end{equation}
From~\eqref{equation} we get
\begin{eqnarray}
  &&\hspace{-6em} \esp\bigl( \e^{-tY} \mid N,\ (W_j)_{j\ge 1}\bigr) = \prod_{1\le
     j\le N} g\bigl( tW_j/b\bigr) \nonumber\\
  &&\hspace{-6em} \esp\bigl( \e^{-tY} \mid N\bigr) = \left( \int g\bigl(
     tx/b\bigr)\,\dif \mu(x)\right)^N\nonumber\\
  &&\hspace{-6em} g(t) = \varphi\Bigl(\int g\bigl(
     tx/b\bigr)\,\dif \mu(x)\Bigr). \label{eqLaplace}
\end{eqnarray}

Let $\displaystyle u(t) = \int g\bigl( tx/b\bigr)\,\dif \mu(x)$. So,
$g(t) = \varphi\bigl(u(t)\bigr)$.

Relations~\eqref{init} become
\begin{equation}\label{init2}
  u(0) = 1 \text{\quad and\quad } u'(0) = -\frac{1}{b}.
\end{equation}
Also $g$ and $u$ are decreasing,
$\displaystyle \lim_{t\to +\infty} g(t) = \proba(Y=0)= \beta$, 
$ \displaystyle \lim_{t\to +\infty} u(t) = \varphi^{-1}(\beta) =
\alpha\beta+1-\alpha$, and  $g(0)=u(0)= 1$.

\subsection*{A particular case}

Now we consider the particular case when
$$
\dif \mu(x) = (1-\alpha )\delta(\dif x) + \alpha (1-\gamma)
b^{\gamma-1}x^{-\gamma}{\large\bf 1}_{(0,b)}(x)\,\dif x,
$$
with $\gamma<1$, $0<\alpha \le 1$, and where~$\delta$ stands for the
unit Dirac mass at~0. If $W$ is distributed according to~$\mu$, the
condition $\esp W = 1$, means $\gamma = 1-\frac1{\alpha b-1}$ and
$\esp W^2 < b$ means $\gamma>1-\frac2{\alpha b-1}$. So the only
constraints on the paramters are
$$\frac1b< \alpha\le 1 \ \text{ and}\quad \gamma = 1-\frac{1}{\alpha b-1}.
$$
\\[1em]
We have
$\displaystyle u(t) = 1-\alpha + \alpha (1-\gamma)b^{\gamma-1}\int_0^b
g\bigl( tx/b\bigr)x^{-\gamma}\dif x$.\\[3pt]
Then
\begin{eqnarray*}
  u'(t) &=& \alpha (1-\gamma)b^{\gamma-1} \int_0^b g'(tx/b)b^{-1}x^{1-\gamma}\dif
            x\\
        &=& \frac{\alpha (1-\gamma)}{t}g(t)-\frac{(1-\gamma)^2b^{\gamma-1}}{t} \int_0^b
            g(tx/b)x^{-\gamma} \dif x\\
        &=& \frac{\alpha (1-\gamma)}{t}g(t) - \frac{1-\gamma}{t}\bigl(u(t)-(1-\alpha)\bigr).
\end{eqnarray*}
We see that $u$ satisfies the following differential equation
\begin{equation}\label{ode}
  (\alpha b-1)u'(t) = \frac{1}{t}\Bigl( \alpha\varphi\bigl(u(t)\bigr)-u(t) + 1 -\alpha\Bigr).
\end{equation}

Let
$\displaystyle \omega(x) = \frac{\alpha b-1}{\alpha \varphi(x)-x + 1 -
  \alpha}+\frac{1}{1-x}$. As
$\varphi(x) = 1 + b(x-1) + \mathrm{O}\bigl((x-1)^2\bigr)$, $\omega$ is
bounded in a neighborhood of~1. Indeed, by continuity, we have
$\displaystyle \omega(1) = -\frac{\esp N(N-1)}{2(\alpha b
  -1)}$. Then~\eqref{ode} rewrites as
\begin{equation}
\left(\omega(u)-\frac{1}{1-u}\right)\dif u = \frac{\dif t}{t}.
\end{equation}

Let $\displaystyle \Omega(x) = \exp \int_1^x
\omega(\tau)\,\dif\tau$. The function
$x \mapsto \alpha \varphi(x)-x + 1 - \alpha$ is convex and vanishes for
$x=1$ and $x=\varphi^{-1}(\beta)$; so it is negative on the interval
$\bigl(\varphi^{-1} (\beta),1\bigr)$. This means that, on this
interval, the derivative of $(1-x)\Omega(x)$ is negative.

It follows that there is a constant~$c$ such that, for $t\ge 1$,
$u(t)$ is the unique solution to equation
$$\bigl(1-u(t)\bigr)\Omega\bigl(u(t)\bigr) = ct, \text{\quad
  with\quad} u(0)=1,$$
in the interval~$\bigl(\varphi^{-1} (\beta),1\bigr)$.
By taking in account the initial conditions~\eqref{init2} we see that
$c=1/b$. Finally, $u$ is implicitly defined by
$$ b(1-u)\Omega(u)=t.$$

\subsection*{Examples}

We give six examples of computations. The fourth one is interesting
because it shows that the mapping $(N,\mu)\longmapsto \nu$ is not
one-to-one.

\subsubsection*{Example 1}

We take  $\varphi(x) = x^{n+1}$ (where ~$n\ge 1$
is an integer) and $\alpha=1$. Then $b=n+1$,
$\gamma=1-\frac1n$, $\beta=0$, and
\begin{eqnarray*}
&&(1-x)\Omega(x) = \frac{1-x^n}{nx^n},\\
&&u(t) = \Bigl(1+\frac{nt}{n+1}\Bigr)^{-1/n},\\
&&g(t) =
\Bigl(1+\frac{nt}{n+1}\Bigr)^{-(n+1)/n}.
\end{eqnarray*}

This means that in this case the variable $Y$ follows the
$\displaystyle\Gamma\biggl( \frac{n+1}{n},\frac{n+1}{n}\biggr)$
distribution, i.e.,
$$\nu(\dif s) =  \Gamma \left(\frac{n+1}{n}\right)^{-1} \left(\frac{n+1}{n}\right)^{\frac{n+1}{n}}s^\frac1n \exp {\displaystyle -\frac{(n+1)s}{n}}\ \dif s. $$
\medskip

This situation has been independently studied by G.~Letac and the
author. It is mentioned in~\cite{G} (page 264), but up to now (\cite{M3},
p.387--388) seems to be the only written trace of this formula.
\medskip

\subsubsection*{Example 2}
This time $\varphi(x) = 1-\rho+\rho x^2$ with $0< \rho\le 1$.\\
Then $b=2\rho\ \text{(therefore $\alpha>\frac1{2\rho}$)},\ \beta = 1-\frac{2\alpha\rho-1}{\rho\alpha^2},\ \gamma = 1- \frac{1}{2 \alpha \rho -1},\\
\omega(x) = -\frac{\alpha \rho}{\alpha \rho x +\alpha \rho -1}$,
\begin{eqnarray*}
  \Omega(x) &=& \frac{2\alpha\rho-1}{\alpha\rho x+\alpha-1},\\
  u(t) &=& \frac{1-\alpha\rho}{\alpha\rho} + \frac{2(2\alpha\rho-1)^2}{ \alpha\rho(4 \alpha\rho +\alpha t-2)},\\
  g(t) &=& \frac{\alpha^2\rho-2\alpha\rho+1}{\alpha^2\rho} + \frac{4(1-\alpha\rho)(2\alpha\rho-1)^2}{\rho\alpha^2(4\alpha\rho+\alpha t-2)}+\frac{4(2 \alpha\rho-1)^{4}}{\rho\alpha^2(4\alpha\rho+\alpha t-2)^2},
\end{eqnarray*}
and
\begin{multline*}
  \nu(\dif s) =
  \frac{(\alpha^2\rho-2\alpha\rho+1)}{\rho\alpha^2}\,\delta(\dif s)\\
  +\frac{4(2\alpha\rho-1)^2\bigl((2\alpha\rho-1)^2s+\alpha(1-\alpha\rho)\bigr)}{\rho
    \,\alpha^{4}}\,{\mathrm e}^{-4 s \rho +\frac{2 s}{\alpha}}\,\dif s.
\end{multline*}\medskip

\subsubsection*{Example 3}
{$\varphi = (1-\rho) x+\rho x^{n+1}$, ($0< \rho\le 1$)
  and $\alpha=1$}.
Then $b=\rho n+1$, $\gamma = (\rho n-1)/\rho n$,
$(1-x)\Omega(x) = (x^{-n}-1)/n,$
\begin{eqnarray*}
u(t) &=& \Bigl(1+\frac{nt}{\rho n+1}\Bigr)^{-1/n},\text{ and}\\
g(t) &=&(1-\rho)\Bigl(1+\frac{nt}{\rho n+1}\Bigr)^{-1/n}
  +\rho\Bigl(1+\frac{nt}{\rho n+1}\Bigr)^{-(n+1)/n}.
\end{eqnarray*}
Finally $\nu$ is a barycenter of Gamma distributions.\medskip

\subsubsection*{Example 4}
This time $\varphi(x) = \displaystyle \frac{1-p}{1-px}$, with
$p>1/2$.\\[3pt]
Then $b = \frac{p}{1-p}\ \text{(so $\alpha> \frac{1-p}{p}$)},\ \gamma
= 1-\frac{1-p}{\alpha p +p -1},\ \omega = -\frac{\alpha
  \,p^{2}}{\left(\alpha p +p x -1\right) \left(1-p \right)}$, and
$$\Omega(x) = \left(\frac{\alpha  p +p x -1}{\alpha  p +p -1}\right)^{-\frac{\alpha  p}{1-p}}.$$

By taking $\alpha = 2(1-p)/p$ the calculation can be pushed
forward. In this condition
\begin{eqnarray*}
  \Omega(x) &=& \left(\frac{1-p}{1+p(x-2)}\right)^2,\\
  u(t) &=& \frac{2p-1}{p}+\frac{2(1-p)}{p(\sqrt{4 t +1}+1)},\\
  g(t) &=& \frac{1}{2}+\frac{1}{2 \sqrt{4 t +1}},\\
  \nu(\dif s) &=& \frac12\delta(\dif s)+\frac{{\mathrm e}^{-\frac{s}{4}}}{4 \sqrt{\pi  s}}\,\dif s.
\end{eqnarray*}
\medskip
                  
\subsubsection*{Example 5}
This time $\varphi = \displaystyle \frac{(1-p)x}{1-px}$, with
$0< p< 1$ and $\alpha > 1-p$. Then
$b = 1+\frac{p}{1-p},\ \gamma = 1-\frac{1-p}{\alpha -1+p}$, and
\begin{eqnarray*}
  \omega(x) &=& -\frac{\alpha  p}{\left(p x +\alpha -1\right) \left(1-p \right)},\\ \Omega(x) &=& \left(\frac{p x +\alpha -1}{p+\alpha -1}\right)^{-\frac{\alpha}{1-p}}
\end{eqnarray*}
Let us end the computation in case when $\alpha = 2(1-p)$. Then we have
\begin{eqnarray*}
  \Omega(x) &=& \left(\frac{1-p}{1+p(x -2)} \right)^{2}\\
  u(t) &=& \frac{2p-1}{p}+\frac{2(1-p)}{p(\sqrt{4 p t +1}+1)},\\
  g(t) &=& \frac{2p-1}{2 p}+\frac{1}{2p\sqrt{4 p t +1}},\\
  \nu(\dif x) &=& \frac{2p-1}{2p}\,\delta(\dif s)+\frac{{\mathrm e}^{-\frac{s}{4 p}}}{4 p^{\frac{3}{2}} \sqrt{\pi  s}}\,\dif s.
\end{eqnarray*}
\medskip

\subsubsection*{Example 6}
This time $\displaystyle \varphi(x) = \frac{(1-p)x^2}{1-px}$. Then
$\displaystyle b=2+\frac{p}{1-p}\text{ and }\alpha
>1-\frac{1}{2-p}$. Also\\
$\gamma =1+\frac{1-p}{\alpha p -2 \alpha -p +1},\ \omega =
-\frac{\alpha}{\left(\left(\alpha +p -\alpha p\right) x +\alpha
    -1\right) \left(1-p \right)}$, and
$$
\Omega(x) = \left(\frac{(\alpha +p -\alpha p) x +\alpha -1}{2\alpha+p -\alpha p-1}\right)^{-\frac{\alpha}{(1-p)(\alpha+p-\alpha p)}}.$$

Now take $\displaystyle \alpha = -\frac{2 p(1-p)}{2 p^{2}-4 p
  +1}$. Due to $1-1/(2-p)< \alpha \le 1$, we have to assume that
$p>1/2$. Then $\gamma = 2(1-p)^2$,
\begin{eqnarray*}
  \Omega(x) &=& \left(\frac{1-p}{1-p(x-2)}\right)^2,\\
  u(t) &=& \frac{2 p -1}{p} + \frac{2(1- p)}{p\left(1+\sqrt{\frac{t}{2-p}+1}\right)},\\
  g(t) &=& \frac{\left(2 p -1\right)^{2}}{2 p^{2}}+\frac{1}{2 p^2\sqrt{\frac{t}{2-p}+1}}-\frac{2(1-p)^{2}}{p^{2} \left(1+\sqrt{\frac{t}{2-p}+1}\right)},
\end{eqnarray*}
and
\begin{multline*}
  \nu(\dif s) = \frac{(2 p -1)^{2}}{2p^2}\,\delta(\dif s )+ \frac{(2 p
    -1) (3-2p)}{2 p^{2}}\,\sqrt{\frac{2-p}{\pi s}}\, {\mathrm
    e}^{(-2+p ) s}\,\dif s \\+\frac{2(2-p ) (1-p
    )^{2}}{p^{2}}\mathrm{erfc}\bigl(\sqrt{(2-p)s}\,\bigr)\,\dif s.
\end{multline*}

\end{document}